\documentclass[a4paper]{article}
\usepackage{amsfonts}
\usepackage[leqno]{amsmath}
\usepackage{amssymb}
\usepackage{amsthm}
\usepackage[utf8]{inputenc} 
\usepackage{graphicx} 
\usepackage[titles,subfigure]{tocloft} 
\usepackage[nottoc]{tocbibind}

\usepackage{booktabs} 
\usepackage{array} 
\usepackage{paralist} 
\usepackage{verbatim} 
\usepackage{subfig} 
\usepackage{soul} 
\usepackage{ngerman} 
\usepackage{textcomp} 
\usepackage{color}
\definecolor{green}{RGB}{40,160,80}
\definecolor{red}{RGB}{192,96,48}
\usepackage{hyperref} 
\hypersetup{
  citecolor  = black 
}
\let\OLDthebibliography\thebibliography
\renewcommand\thebibliography[1]{
  \OLDthebibliography{#1}
  \setlength{\parskip}{0pt}
  \setlength{\itemsep}{10pt plus 0.3ex}
}

\usepackage{fancyhdr} 
\DeclareMathOperator{\ab}{ab}

\DeclareMathOperator{\End}{End}
\DeclareMathOperator{\Gal}{Gal}
\DeclareMathOperator{\GL}{GL}
\DeclareMathOperator{\hgt}{ht}
\DeclareMathOperator{\Hom}{Hom}

\DeclareMathOperator{\mSpec}{m-Spec}

\DeclareMathOperator{\nr}{nr}
\DeclareMathOperator{\Spec}{Spec}

\DeclareMathOperator{\tr}{tr}
\DeclareMathOperator{\univ}{univ}

\newcommand{\C}{\mathbb{C}}

\newcommand{\m}{\mathfrak{m}}
\newcommand{\M}{\mathrm{M}}
\newcommand{\N}{\mathbb{N}}
\newcommand{\Oc}{\mathcal{O}}

\newcommand{\p}{\mathfrak{p}}
\newcommand{\q}{\mathfrak{q}}
\newcommand{\Q}{\mathbb{Q}}

\newcommand{\x}{\tilde{x}}
\newcommand{\y}{\tilde{y}}
\newcommand{\z}{\tilde{z}}
\newcommand{\Z}{\mathbb{Z}}
\newcommand{\1}{\mathbf{1}}

\newtheorem{prop}{Proposition}[section]

\newtheorem{lem}[prop]{Lemma}
\newtheorem{thm}[prop]{Theorem}
\newtheorem{satz}[prop]{Satz}

\newtheorem{nota}[prop]{Notation}
\newtheorem{kor}[prop]{Korollar}
\newenvironment{bew}{\begin{proof}[Beweis]}{\end{proof}}
\theoremstyle{definition}
\newtheorem{defi}[prop]{Definition}
\theoremstyle{remark}
\newtheorem{bem}[prop]{Bemerkung}
\newtheorem{beisp}[prop]{Beispiel}

\newcommand{\OO}{\Oc}
\newcommand{\rhobar}{\bar{\rho}}
\newcommand{\Eins}{\mathbf{1}}
\newcommand{\mm}{\mathfrak m}

\title{Irreduzible Komponenten von $2$-adischen Deformationsräumen}

\author{Maurice Babnik}

\begin{document}
\maketitle
    
\abstract{ We prove that the irreducible components of the space of framed deformations of a  $2$-dimensional mod $2$ representation with scalar semi-simplification of the absolute Galois group of $\Q_2$ are in natural bijection with those of its determinant, confirming a conjecture of B\"ockle--Juschka. }
\tableofcontents


\section{Introduction}

Let $\Q_2$ be the field of $2$-adic numbers and let $G_{\Q_2}$ be its absolute Galois group. Let $L$ be a finite field extension of $\Q_2$ with the
ring of integers $\OO$ and residue field $k$.  

Let $\rhobar: G_{\Q_2}\rightarrow \GL_2(k)$ be a continuous representation of the form 
$\rhobar \cong \bigl(\begin{smallmatrix} 1 & \ast \\ 0 & 1 \end{smallmatrix}\bigr)$.   Let $\Eins$ be a one dimensional $k$-vector space on which $G_{\Q_2}$ acts trivially, and let $D_{\Eins}$ be the deformation functor of $\Eins$, and let $D^{\Box}$ 
be the  framed deformation functor of $\rhobar$, so that for each local artinian $\OO$-algebra $(A, \mm_A)$ with residue field $k$, $D_{\Eins}(A)$ is the set of continuous group 
homomorphisms from $G_{\Q_2}$ to $1+\mm_A$ and $D^{\Box}(A)$ is the set of continuous group 
homomorphisms from $\rho_A:G_{\Q_2}\rightarrow \GL_2(A)$, such by reducing the matrix entries of $\rho_A$ modulo $\mm_A$ we obtain $\rhobar$. 

 These functors are represented by complete local noetherian $\OO$-algebras $R_{\Eins}$ and $R^{\Box}$ respectively.
Mapping a framed deformation of $\rhobar$ to its determinant induces a natural transformation $D^{\Box}\rightarrow D_{\Eins}$, and hence a homomorphism of $\OO$-algebras  $d: R_{\Eins}\rightarrow R^{\Box}$. 
 
 \begin{thm}\label{vienas} The map $d: R_{\Eins}\rightarrow R^{\Box}$ induces a bijection between the irreducible components of $\Spec R^{\Box}$ and $\Spec R_{\Eins}$. In particular, $\Spec R^{\square}$ has two irreducible components. 
 \end{thm}

Our main result answers affirmatively a question of B\"ockle-Juschka \cite{BJ} in this case. 

The proof closely follows the strategy employed
by Colmez--Dospinescu--Pa\v{s}k\={u}nas in \cite{CDP}, where the case $\rhobar \cong \bigl(\begin{smallmatrix} 1 & 0 \\ 0 & 1 \end{smallmatrix}\bigr)$, so that $\rhobar$ is split, is considered. We will briefly recall the strategy. We first show that any framed deformation $\rho_A: G_{\Q_2}\rightarrow \GL_2(A)$ factors through the maximal pro-$2$ quotient $G_{\Q_2}(2)$ of $G_{\Q_2}$. This group is known to be topologically generated by three generators, 
which satisfy one relation. Using this we present the ring $R^{\square}$ in Satz \ref{rsquareisom} as a quotient of the ring of formal power series 
over $\OO$ in $12$ variables by $4$ relations. We call this presentation $S$. We show in Lemma \ref{vollstdurchschn} that $S$ is complete intersection, which implies that $S[1/2]$ is Cohen--Macaulay. In Proposition \ref{singort} we bound the dimension of the singular locus in $S[1/2]$.
Using Serre's criterion for normality we deduce that $S[1/2]$ is a product of normal domains. This part of the proof is essentially the same as in 
\cite{CDP}, albeit our regular sequence in Lemma \ref{vollstdurchschn} is different to the one considered in \cite{CDP}. We then show that the 
closed points of $\Spec S[1/2]$, which are expected to lie on the same irreducible component, can be connected by a sequence of $p$-adic discs.
If this is possible then the theory developed in \cite{CDP} implies that such points lie on the same component. 
The main difficulty is to actually produce such disks.  

Overcoming this problem is the most original part of the paper. The construction of disks in  \cite{CDP} does not seem 
to carry over directly, when $\rhobar$ is non-split. If $\rhobar$ is split  then it is easy to produce lifts of $\rhobar$ to characteristic $0$ by writing down diagonal matrices for the generators and this makes it  easier to write down explicit formulas for the disks connecting such points; this is exploited in \cite{CDP}.  If $\rhobar$ is non-split then such diagonal lifts do not exist. We first show that every irreducible component of $S[1/2]$ contains a point, where the matrices for the 
generators are upper triangular, by constructing an explicit  regular sequence in $S$, and looking at its zero locus, see
Lemma \ref{schnitt}. We then show directly that the largest quotient of $S[1/2]$, where the all the matrices for the generators are upper-triangular has $4$ irreducible components, see Satz \ref{4irredkomp}. This allows us to produce
families of very special points, which meet every irreducible component of $S[1/2]$, see Korollar \ref{punkte1}, Korollar \ref{punkte2}.
We then write down explicit formulas for the disks connecting the points, which are expected to lie on the same
irreducible component,  see \S \ref{disks}.

We state the results for $\rhobar$ of the form $\bigl(\begin{smallmatrix} 1 & \ast \\ 0 & 1 \end{smallmatrix}\bigr)$, but our results apply to any representation of the form $\bigl(\begin{smallmatrix} \chi & \ast \\ 0 & \chi \end{smallmatrix}\bigr)$, where $\chi: G_{\Q_2}\rightarrow k^{\times}$ is any continuous character, as can easily be seen by twisting with 
a lift of $\chi$. 

This paper is a shortened version of my Master thesis at Universit\"at Duisburg--Essen written under the direction of 
Prof. Dr. V. Pa\v{s}k\={u}nas.

\textbf{Danksagung.} Mein besonderer Dank für die Betreuung und Begleitung meiner Masterarbeit gilt Prof. Dr. Vytautas Pa\v{s}k\=unas, dem Zweitgutachter Prof. Dr. Georg Hein, ebenso meiner ganzen Familie, meinem Bruder Etienne Jerôme Babnik und \mbox{Frau Cornelia Steinbock.}

\newpage

\section{Der Ring \texorpdfstring{$S$}{S}}

\begin{defi}
Nun setzen wir
$$A \mathrel{\mathop:}= \Oc[[X, Y, Z]] \mathrel{\mathop:}= \Oc[[x_{11}, x_{12}, x_{21}, x_{22}, y_{11}, y_{12}, y_{21}, y_{22}, z_{11}, z_{12}, z_{21}, z_{22}]]$$
(der Ring der formalen Potenzreihen über $\Oc$ in 12 Variablen); dieser Ring ist ein $13$-dimensionaler vollständiger regulärer lokaler Ring zum maximalen Ideal $\m_A \mathrel{\mathop:}= (\varpi, x_{11}, \ldots , z_{22})$, und somit insbesondere noethersch. Außerdem wählen wir $\lambda, \mu, \kappa \in \Oc$. Nun definieren wir Matrizen $\tilde X, \tilde Y, \tilde Z \in \M_{2 \times 2}(A)$ wie folgt:
\begin{align*}
\tilde X &\mathrel{\mathop:}= \begin{pmatrix} 1 + x_{11} & \lambda + x_{12} \\ x_{21} & 1 + x_{22} \end{pmatrix} =\mathrel{\mathop:} \begin{pmatrix} \x_{11} & \x_{12} \\ \x_{21} & \x_{22} \end{pmatrix}, \\ \tilde Y &\mathrel{\mathop:}= \begin{pmatrix} 1 + y_{11} & \mu + y_{12} \\ y_{21} & 1 + y_{22} \end{pmatrix} =\mathrel{\mathop:} \begin{pmatrix} \y_{11} & \y_{12} \\ \y_{21} & \y_{22} \end{pmatrix}, \\ \text{und } \tilde Z &\mathrel{\mathop:}= \begin{pmatrix} 1 + z_{11} & \kappa + z_{12} \\ z_{21} & 1 + z_{22} \end{pmatrix} =\mathrel{\mathop:} \begin{pmatrix} \z_{11} & \z_{12} \\ \z_{21} & \z_{22} \end{pmatrix}.
\end{align*}
Damit erhalten wir:
\begin{align*}
\tilde X \bmod \m_A &= \begin{pmatrix} 1 & \bar\lambda \\ 0 & 1 \end{pmatrix}, \ \tilde Y \bmod \m_A = \begin{pmatrix} 1 & \bar\mu \\ 0 & 1 \end{pmatrix}, \\ \text{und } \tilde Z \bmod \m_A &= \begin{pmatrix} 1 & \bar\kappa \\ 0 & 1 \end{pmatrix} \in \M_{2 \times 2}(k).
\end{align*}
Wir setzen außerdem
$$\tilde X^2 \tilde Y^4 [\tilde Y, \tilde Z] =\mathrel{\mathop:} \begin{pmatrix} 1 + f_{11} & f_{12} \\ f_{21} & 1 + f_{22} \end{pmatrix}, \text{ wobei } [A, B] \mathrel{\mathop:}= ABA^{-1}B^{-1}.$$
Sei nun schließlich $S \mathrel{\mathop:}= S_{\lambda, \mu, \kappa} \mathrel{\mathop:}= A / (f_{11}, f_{12}, f_{21}, f_{22})$; dann ist auch $S$ (als Quotient eines vollständigen Rings) vollständig, lokal und noethersch.
\end{defi}

\begin{bem}
In der ursprünglichen Arbeit \cite{CDP} galt $\lambda = \mu = \kappa = 0$. 
\end{bem}

\begin{lem}
Seien $\lambda_1, \mu_1, \kappa_1, \lambda_2, \mu_2, \kappa_2 \in \Oc$, so dass $\lambda_1 \equiv \lambda_2 \pmod {(\varpi)}$, $\mu_1 \equiv \mu_2 \pmod {(\varpi)}$ und $\kappa_1 \equiv \kappa_2 \pmod {(\varpi)}$. Dann gibt es einen Isomorphismus $S_{\lambda_1, \mu_1, \kappa_1} \cong S_{\lambda_2, \mu_2, \kappa_2}$. Insbesondere können wir o.~B.~d.~A. annehmen, dass $\lambda, \mu, \kappa$ entweder Einheiten in $\Oc$ oder $= 0$ sind.
\end{lem}

\begin{bew}
Da $\lambda_2 - \lambda_1, \mu_2 - \mu_1, \kappa_2 - \kappa_1 \in (\varpi)$, induziert die Abbildung
$$x_{12} \mapsto x_{12} + (\lambda_2 - \lambda_1), \quad y_{12} \mapsto y_{12} + (\mu_2 - \mu_1), \quad z_{12} \mapsto z_{12} + (\kappa_2 - \kappa_1)$$
einen Automorphismus von $A$ (wie in Kapitel 2 definiert), der $x_{12}$, $y_{12}$ und $z_{12}$ auf $x_{12} + (\lambda_2 - \lambda_1)$, $y_{12} + (\mu_2 - \mu_1)$ und $z_{12} + (\kappa_2 - \kappa_1)$, also $\lambda_1 + x_{12}$, $\mu_1 + y_{12}$ und $\kappa_1 + z_{12}$ auf $\lambda_2 + x_{12}$, $\mu_2 + y_{12}$ und $\kappa_2 + z_{12}$ abbildet, und alle formalen Potenzreihen, in denen $x_{12}$, $y_{12}$ und $z_{12}$ nicht vorkommen, fest lässt. Dieser Automorphismus induziert dann einen Isomorphismus $S_{\lambda_1, \mu_1, \kappa_1} \cong S_{\lambda_2, \mu_2, \kappa_2}$.
\end{bew}

\begin{defi}
Sei außerdem $G_{\Q_2}(2)$ der maximale pro-$2$-Quotient von $G_{\Q_2}$.  Wir fixieren einen algebraischen Abschluss $\overline{\Q_2}$ von $\Q_2$. Nun definieren wir $\overline{\Q_2}(2)$ als die Vereinigung der Körper $K$ mit $\Q_2 \subseteq K \subseteq \overline{\Q_2}$, so dass $K$ galoissch über $\Q_2$ ist, und $[K : \Q_2]$ eine Zweierpotenz ist. Dann ist $\overline{\Q_2}(2)$ galoissch über $\Q_2$, und wir haben $G_{\Q_2}(2) = \Gal(\overline{\Q_2}(2)/\Q_2)$.
\end{defi}

\begin{satz}
\label{rsquareisom}
Es gibt einen Isomorphismus $R^\square \cong S$ von $\Oc$-Algebren.
\end{satz}

\begin{bew}
Nach \cite{Ser} wird die Gruppe $G_{\Q_2}(2)$ durch $3$ Erzeuger $x$, $y$ und $z$ mit der Relation $x^2 y^4 [y, z] = 1$ topologisch erzeugt. Dann wählen wir $\lambda, \mu, \kappa \in \Oc$ so, dass $\bar\rho(x) = \begin{pmatrix} 1 & \bar\lambda \\ 0 & 1 \end{pmatrix}$, $\bar\rho(y) = \begin{pmatrix} 1 & \bar\mu \\ 0 & 1 \end{pmatrix}$ und $\bar\rho(z) = \begin{pmatrix} 1 & \bar\kappa \\ 0 & 1 \end{pmatrix}$ gelten.

Da $\tilde X$, $\tilde Y$ und $\tilde Z$ modulo $\m_S$ sich zu $\begin{pmatrix} 1 & \bar\lambda \\ 0 & 1 \end{pmatrix}$, $\begin{pmatrix} 1 & \bar\mu \\ 0 & 1 \end{pmatrix}$ und $\begin{pmatrix} 1 & \bar\kappa \\ 0 & 1 \end{pmatrix}$ reduzieren, ist damit durch $x \mapsto \tilde X$, $y \mapsto \tilde Y$, $z \mapsto \tilde Z$ eine eindeutige stetige Darstellung ${\rho_S : G_{\Q_2}(2) \rightarrow \GL_2(S)}$ definiert. Wir betrachten $\rho_S$ als eine gerahmte Deformation von $\bar\rho$ zu $(S, \m_S)$. Dies liefert einen Homomorphismus $\varphi: R^\square \rightarrow S$ von $\Oc$-Algebren.

Sei $(A, \m_A)$ eine lokale artinsche $\Oc$-Algebra mit Restklassenkörper $k$. Die Menge $D^\square(A)$ steht dann in Bijektion zur Menge der stetigen Gruppenhomomorphismen $\rho : G_{\Q_2} \rightarrow \GL_2(A)$ mit $\rho_A \mod \m_A = \bar\rho$. Da die Untergruppe der Matrizen in $\GL_2(A)$, die sich modulo $\m_A$ zu $\begin{pmatrix} 1 & * \\ 0 & 1 \end{pmatrix}$ reduzieren, Ordnung $2^n$ für ein $n \in \N_0$ hat, faktorisiert jedes solches $\rho$ durch den maximalen pro-$2$ Quotient $G_{\Q_2}(2)$ von $G_{\Q_2}$. Damit induziert die Abbildung, die jedes $\rho : G_{\Q_2}(2) \rightarrow \GL_2(A)$ auf
$$\left(\rho(x) - \begin{pmatrix} 1 & \lambda \\ 0 & 1 \end{pmatrix}, \rho(y) - \begin{pmatrix} 1 & \mu \\ 0 & 1 \end{pmatrix}, \rho(z) - \begin{pmatrix} 1 & \kappa \\ 0 & 1 \end{pmatrix}\right)$$
abbildet, eine Bijektion zwischen der Menge der solchen $\rho$ und der Menge der Tripel $(X_A, Y_A, Z_A) \in \M_2(\m_A)^3$, welche $\tilde X_A^2 \tilde Y_A^4 [\tilde Y_A, \tilde Z_A] = 1$ erfüllen, wobei $\tilde X_A = \begin{pmatrix} 1 & \lambda \\ 0 & 1 \end{pmatrix} + X_A$, $\tilde Y_A = \begin{pmatrix} 1 & \mu \\ 0 & 1 \end{pmatrix} + Y_A$ und $\tilde Z_A = \begin{pmatrix} 1 & \kappa \\ 0 & 1 \end{pmatrix} + Z_A$. Diese stehen wiederum in Bijektion zur Menge $\Hom_\Oc(S, A)$. Damit ist $\varphi: R^\square \rightarrow S$ ein Isomorphismus.
\end{bew}

%
%

Sei $\delta \mathrel{\mathop:}= \det \tilde X (\det \tilde Y)^2 \in S$, dann liefert die Relation in $S$, dass $\delta^2 = 1$. Damit erwartet man, dass die irreduziblen Komponenten durch $\delta = 1$ und $\delta = -1$ gegeben sind. Ist $\p \in \Spec S$, dann ist $S / \p$ ein Integritätsbereich. Die Relation $(\delta - 1) (\delta + 1) = \delta^2 - 1 = 0$ impliziert, dass $\delta = 1$ oder $\delta = -1$ in $S / \p$ gelten muss.

\begin{nota}
Wir setzen $S^+ \mathrel{\mathop:}= S / (\delta + 1)$ und $S^- \mathrel{\mathop:}= S / (\delta - 1)$. Außerdem schreiben wir $S^\pm$, wenn einer der beiden Ringe $S^+$ oder $S^-$ gemeint ist.
\end{nota}

%
%
%

\begin{lem}
\label{kommut}
Seien
$$A \mathrel{\mathop:}= \begin{pmatrix} a_{11} & a_{12} \\ a_{21} & a_{22} \end{pmatrix}, \text{ und } B \mathrel{\mathop:}= \begin{pmatrix} b_{11} & b_{12} \\ b_{21} & b_{22} \end{pmatrix}$$
zwei Matrizen mit Elementen in einem kommutativen Ring $R$. Dann gilt genau dann $AB = BA$, wenn die Determinanten aller $2 \times 2$-Untermatrizen in $\begin{pmatrix} a_{11} - a_{22} & a_{12} & a_{21} \\ b_{11} - b_{22} & b_{12} & b_{21} \end{pmatrix}$ verschwinden.
\end{lem}

\begin{bew}
Wir haben
$$AB = \begin{pmatrix} a_{11} & a_{12} \\ a_{21} & a_{22} \end{pmatrix} \begin{pmatrix} b_{11} & b_{12} \\ b_{21} & b_{22} \end{pmatrix} = \begin{pmatrix} a_{11} b_{11} + a_{12} b_{21} & a_{11} b_{12} + a_{12} b_{22} \\ a_{21} b_{11} + a_{22} b_{21} & a_{21} b_{12} + a_{22} b_{22} \end{pmatrix}$$
und
$$BA = \begin{pmatrix} b_{11} & b_{12} \\ b_{21} & b_{22} \end{pmatrix} \begin{pmatrix} a_{11} & a_{12} \\ a_{21} & a_{22} \end{pmatrix} = \begin{pmatrix} a_{11} b_{11} + a_{21} b_{12} & a_{12} b_{11} + a_{22} b_{12} \\ a_{11} b_{21} + a_{21} b_{22} & a_{12} b_{21} + a_{22} b_{22} \end{pmatrix},$$
also
$$AB - BA = \begin{pmatrix} a_{12} b_{21} - a_{21} b_{12} & (a_{11} - a_{22}) b_{12} - a_{12} (b_{11} - b_{22}) \\ a_{21} (b_{11} - b_{22}) - (a_{11} - a_{22}) b_{21} & a_{21} b_{12} - a_{12} b_{21} \end{pmatrix};$$
somit gilt genau dann $AB = BA$, wenn die Determinanten der drei Matrizen $\begin{pmatrix} a_{12} & a_{21} \\ b_{12} & b_{21} \end{pmatrix}$, $\begin{pmatrix} a_{11} - a_{22} & a_{12} \\ b_{11} - b_{22} & b_{12} \end{pmatrix}$ und $\begin{pmatrix} a_{11} - a_{22} & a_{21} \\ b_{11} - b_{22} & b_{21} \end{pmatrix}$ verschwinden.
\end{bew}

\begin{lem}
\label{specspm12}
$\Spec S^\pm[1/2]$ ist eine abgeschloffene, d.~h. abgeschlossene und offene Teilmenge von $\Spec S[1/2]$ und eine (endliche) Vereinigung irreduzibler Komponenten von $\Spec S[1/2]$.
\end{lem}

\begin{bew}
Wegen $\delta^2 = 1$ haben wir $(\frac{1 \pm \delta}2)^2 = \frac{1 \pm 2 \delta + \delta^2}4 = \frac{2 \pm 2 \delta}4 = \frac{1 \pm \delta}2$, d.~h. $\frac{1 \pm \delta}2$ ist ein idempotentes Element in $S[1/2]$. Dessen Nullort ist, da $2$ eine Einheit in $S[1/2]$ ist, gleich dem Nullort von $1 \pm \delta$, nämlich gleich
$$\Spec S[1/2] / (1 \pm \delta) = \Spec (S / (1 \pm \delta))[1/2] = \Spec S^\pm[1/2].$$
Da die beiden idempotenten Elemente $\frac{1 \pm \delta}2$ sich zu $1$ addieren, folgt somit $S[1/2] \cong S^+[1/2] \times S^-[1/2]$, und daraus folgt die Behauptung.
\end{bew}

\begin{lem}
\label{vollstdurchschn}
$S$ ist ein vollständiger Durchschnitt von Dimension $9$, und die Elemente $\varpi$, $x_{11}$, $x_{12}$, $x_{21} + \tr \tilde X$, $y_{21}$, $z_{21}$ bilden eine reguläre Folge in $S$.
\end{lem}

\begin{bew}
Da $A$ ein $13$-dimensionaler regulärer lokaler Ring ist, folgt mit \cite[Thm. 17.8]{Mat}, dass $A$ \textsc{Cohen-Macaulay} ist; außerdem ist $A$ noethersch. Falls wir zeigen können, dass der Quotient von $A$ durch das von den zehn Elementen
$$f_{11}, f_{12}, f_{21}, f_{22}, \varpi, x_{11}, x_{12}, x_{21} + \tr \tilde X, y_{21}, z_{21},$$
die alle in $\m$ liegen, erzeugte Ideal höchstens $3$-dimensional ist, folgt mit \cite[Thm. 17.4 (i)]{Mat}, dass das genannte Ideal Höhe von mindestens 10 hat, mit \cite[Thm. 13.5]{Mat}, dass das genannte Ideal Höhe von genau $10$ hat, mit \cite[Thm. 17.4 (iii)]{Mat}, dass diese Elemente eine reguläre Folge in $A$ bilden und dass das von den ersten vier Elementen erzeugte Ideal Höhe $4$ hat, mit $S = A / (f_{11}, f_{12}, f_{21}, f_{22})$, dass $S$ $9$-dimensional ist, mit der Definition der regulären Folge, dass die Elemente $\varpi$, $x_{11}$, $x_{12}$, $x_{21} + \tr \tilde X$, $y_{21}$, $z_{21}$ eine reguläre Folge in $S$ bilden, und mit \cite[Thm. 21.2 (ii)]{Mat}, dass $S$ ein vollständiger Durchschnitt ist, also die Behauptung.

Es reicht also zu zeigen, dass der genannte Quotient, der isomorph zum Ring $B \mathrel{\mathop:}= S / (\varpi, x_{11}, x_{12}, x_{21} + \tr \tilde X, y_{21}, z_{21})$ ist, höchstens $3$-dimensional ist, weil $S$ der Quotient des $13$-dimensionalen lokalen regulären Rings $A$ durch die ersten $4$ Elemente ist, und wir durch die weiteren $6$ Elemente, die alle im maximalen Ideal liegen, teilen.

Modulo $(\varpi, x_{11}, x_{12}, x_{21} + \tr \tilde X, y_{21}, z_{21})$ sind $\tilde Y$, $\tilde Z$, und somit auch $\tilde Y^4 [\tilde Y, \tilde Z]$, obere Dreiecksmatrizen. Damit impliziert die Relation $\tilde X^2 \tilde Y^4 [\tilde Y, \tilde Z] = 1$, dass auch $\tilde X^2$ modulo $(\varpi, x_{11}, x_{12}, x_{21} + \tr \tilde X, y_{21}, z_{21})$ eine obere Dreiecksmatrix ist, d.~h. wir haben $x_{21}^2 = x_{21} \tr \tilde X = 0$ in $B$.

Sei nun $\p$ ein Primideal von $B$, welches wir als ein Primideal von
$$k[[X, Y, Z]] \mathrel{\mathop:}= k[[x_{11}, x_{12}, x_{21}, x_{22}, y_{11}, y_{12}, y_{21}, y_{22}, z_{11}, z_{12}, z_{21}, z_{22}]],$$
welches
$$({\tilde X^2 \tilde Y^4 [\tilde Y, \tilde Z] = 1}, x_{11}, x_{12}, {x_{21} + \tr \tilde X}, y_{21}, z_{21})$$
enthält, auffassen. Von nun an rechnen wir modulo $\p$. Da $x_{21}^2 = 0$ in $A$, folgt $\tr \tilde X = x_{21} = 0$, und somit $x_{22} = \tr \tilde X + x_{11} = 0$. Damit haben wir $\tilde X = \begin{pmatrix} 1 & \lambda \\ 0 & 1 \end{pmatrix}$ und $\tilde X^2 = 1$. Also vereinfacht sich die Relation $\tilde X^2 \tilde Y^4 [\tilde Y, \tilde Z] = 1$ zu $\tilde Y^4 [\tilde Y, \tilde Z] = 1$. Betrachten wir die Hauptdiagonalelemente, so erhalten wir mit dem Frobeniushomomorphismus $\y_{11}^4 = 1 + y_{11}^4 = \y_{22}^4 = 1 + y_{22}^4 = 1$, also $y_{11}^4 = y_{22}^4 = 0$ und somit $y_{11} = y_{22} = 0$. Damit haben wir $\tilde Y = \begin{pmatrix} 1 & \y_{12} \\ 0 & 1 \end{pmatrix}$ und $\tilde Y^4 = 1$. Somit vereinfacht sich die Relation zu $[\tilde Y, \tilde Z] = 1$ bzw. $\tilde Y \tilde Z = \tilde Z \tilde Y$. Dies impliziert, dass die Determinanten aller $2 \times 2$-Matrizen in $\begin{pmatrix} 0 & \y_{12} & 0 \\ \tr \tilde Z & \z_{12} & 0 \end{pmatrix}$ nach Lemma \ref{kommut} verschwinden, d.~h. wir haben $\y_{12} \tr \tilde Z = 0$. Dies impliziert, dass wir entweder $y_{12} = 0$ oder $\tr \tilde Z = z_{11} + z_{22} = 0$ haben; falls $\mu \neq 0$, tritt stets letzterer Fall ein. Somit schließen wir, dass die Surjektion $k[[X, Y, Z]] \twoheadrightarrow A / \p$ entweder durch $k[[z_{11}, z_{12}, z_{22}]] \twoheadrightarrow B / \p$ oder aber durch $k[[y_{12}, z_{11}, z_{12}]] \twoheadrightarrow B / \p$ faktorisiert. In beiden Fällen erhalten wir $\dim B / \p \leq 3$. Da wir $\p \in \Spec B$ beliebig gewählt haben, folgt $\dim B \leq 3$ und somit die Behauptung.
\end{bew}

\begin{defi}
Sei $\rho^{\univ} : G_{\Q_2} \rightarrow \GL_2(R^\square)$ die universelle gerahmte Deformation zu $\bar\rho$. Ist $x$ ein beliebiger abgeschlossener Punkt von $\Spec R^\square[1/2]$, dann ist sein Restklassenkörper $\kappa(x)$ eine endliche Erweiterung von $L$. Mit $\rho^{\univ}_x : G_{\Q_2} \rightarrow \GL_2(\kappa(x))$ bezeichnen wir die Darstellung, die durch Spezialisierung der universellen Darstellung zu $x$ entsteht.
\end{defi}


\begin{lem}
\label{singpunkt}
Ist $x$ ein abgeschlossener singulärer Punkt von $\Spec S[1/2]$, dann gibt es einen Charakter $\delta : G_{\Q_2} \rightarrow \Oc_{\kappa(x)}^*$ und eine exakte Sequenz
$$0 \rightarrow \delta \rightarrow \rho^{\univ}_x \rightarrow \delta \varepsilon \rightarrow 0$$
von $\kappa(x)[G_{\Q_2}]$-Modulen, wobei $\varepsilon : G_{\Q_2} \rightarrow \Z_2^\times$ der $2$-adische zyklotomische Charakter ist.
\end{lem}

\begin{bew}
Der Beweis geht analog zum Beweis von \cite[Lemma 4.1]{CDP}.
\end{bew}

\begin{lem}
$S[1/2]$ ist exzellent.
\end{lem}

\begin{bew}
Da $A$ ein vollständiger noetherscher lokaler Ring ist, ist $A$ exzellent; da $S$ ein Quotient von $A$ ist, ist $S$ exzellent; und da $S[1/2]$ eine Lokalisierung von $S$ ist, ist $S[1/2]$ exzellent (jeweils nach \cite[Seite 260]{Mat}).
\end{bew}

\begin{prop}
\label{singort}
Der singuläre Ort von $\Spec S[1/2]$ hat Dimension $\leq 6$.
\end{prop}

\begin{bew}
(Nach \cite[Prop. 4.2]{CDP}:)

Da $S[1/2]$ exzellent ist, folgt mit (3) in der Definition auf \cite[Seite 260]{Mat}, dass der singuläre Ort von $\Spec S[1/2]$ abgeschlossen ist, d.h. er kann als $V(I) \cong \Spec S[1/2]/I$ für ein gewisses Ideal $I$ von $A$ geschrieben werden. Da $S[1/2]$ Jacobson ist, impliziert dies, 
dass auch der singuläre Ort Jacobson ist.

Nun folgt mit Lemma \ref{singpunkt}, dass alle singulären abgeschlossenen Punkte von $S[1/2]$ in $V(I)$ enthalten sind, wobei $I$ diesmal das Ideal von $S$ ist, welches durch die Elemente
\begin{equation}
\label{elemente}
(\tr \rho^{\univ}(g))^2 - (\varepsilon(g) + 1)^2 \varepsilon(g)^{-1} \det \rho^{\univ}(g),
\end{equation}
wobei $g$ über $G_{\Q_2}$ variiert, erzeugt wird. Da $S[1/2]/I$ Jacobson ist, ist somit der singuläre Ort auch in $V(I)$ enthalten. Damit genügt es zu zeigen, dass $\dim S/I \leq 7$, weil \cite[Lemma 2.3]{CDP} impliziert, dass dann $\dim (S/I)[1/2] \leq 6$ erfüllt ist.

Sei nun $J \mathrel{\mathop:}= \sqrt{(\varpi, I)}$ und sei $\tilde\rho : G_{\Q_2} \rightarrow GL_2(S/J)$ die Darstellung, die man durch Reduktion der Einträge von $\rho^{\univ}$ modulo $J$ erhält. Da $S/I$ ein lokaler noetherscher Ring ist, folgt mit $\varpi \in \m$, \cite[Thm. 13.6 (ii)]{Mat} und \cite[Seite 3]{Mat}, dass
\begin{displaymath}
\begin{split}
\dim S/I \leq \dim (S/I)/(\varpi) + 1 = \dim S/(\varpi, I) + 1 &= \dim S/\sqrt{(\varpi, I)} + 1 \\ &= \dim S/J + 1.
\end{split}
\end{displaymath}
Somit reicht es aus, $\dim S/J$ durch $6$ nach oben zu beschränken.

Da $\varepsilon(g) \equiv 1 \pmod \varpi$, folgt nun aus (\ref{elemente}), dass $(\tr \rho^{\univ}(g))^2 \equiv 0 \pmod {(\varpi, I)}$, und somit $\tr \tilde\rho(g) = 0$ für alle $g \in G_{\Q_2}$ gilt. Damit faktorisiert die Surjektion $S \twoheadrightarrow S/J$ durch
\begin{equation}
\label{BSJ}
B \mathrel{\mathop:}= \frac{k[[X, Y, Z]]}{(\det \tilde X - \det \tilde Y^{-2}, \tr \tilde X, \tr \tilde Y, \tr \tilde Z, \tr \tilde X \tilde Y, \tr \tilde X \tilde Z, \tr \tilde Y \tilde Z)} \twoheadrightarrow S/J.
\end{equation}
Unter Verwendung der Notation $\x_{12} = \lambda + x_{12}$, $\y_{12} = \lambda + y_{12}$, $\z_{12} = \kappa + z_{12}$ bemerken wir, dass wenn $\tr \tilde Y = \tr \tilde Z = 0$ gilt, daraus $\tr \tilde Y \tilde Z = \y_{12} z_{21} + y_{21} \z_{12}$ folgt, weil wir in Charakteristik $2$ sind. Sei $I'$ das Ideal von $k[x_{12}, x_{21}, y_{12}, y_{21}, z_{12}, z_{21}]$, welches durch alle $2 \times 2$-Minoren der Matrix
$$\begin{pmatrix} \x_{12} & \y_{12} & \z_{12} \\ x_{21} & y_{21} & z_{21} \end{pmatrix}$$
erzeugt wird. Setzen wir $x = \x_{12}$, $y = \y_{12}$ und $z = \z_{12}$, so erhalten wir einen Isomorphismus $k[x_{12}, x_{21}, y_{12}, y_{21}, z_{12}, z_{21}] / I' \cong k[x, x_{21}, y, y_{21}, z, z_{21}] / I''$, wobei $I''$ das Ideal von $k[x, x_{21}, y, y_{21}, z, z_{21}]$ ist, welches durch alle $2 \times 2$-Minoren der Matrix
$$\begin{pmatrix} x & y & z \\ x_{21} & y_{21} & z_{21} \end{pmatrix}$$
erzeugt wird. Nun folgt mit \cite[Proposition 1.1]{BV}, dass $I''$ und somit auch $I'$ eine irreduzible Varietät von Dimension $4$ definiert. Dies impliziert, dass
$$A \mathrel{\mathop:}= \frac{k[[x_{12}, x_{21}, y_{12}, y_{21}, z_{12}, z_{21}]]}{(\x_{12} y_{21} + x_{21} \y_{12}, \x_{12} z_{21} + x_{21} \z_{12}, \y_{12} z_{21} + y_{21} \z_{12})}$$
$4$-dimensional ist. Die Relation $\det \tilde X - \det \tilde Y^{-2}$ impliziert, dass $B$ endlich über $A[[y_{11}, z_{11}]]$ ist, und somit $\dim B = 6$ gilt. Und schließlich folgt aus (\ref{BSJ}), dass $\dim S/J \leq 6$ gilt.
\end{bew}

\begin{lem}
$S[1/2]$ ist \textsc{Cohen-Macaulay}.
\end{lem}

\begin{bew}
Da $S$ ein vollständiger Durchschnitt ist, ist $S$ nach \cite[Seite 171]{Mat} \textsc{Cohen-Macaulay}; und da $S[1/2]$ eine Lokalisierung von $S$ ist, ist $S[1/2]$ nach \cite[Seite 136]{Mat} \textsc{Cohen-Macaulay}.
\end{bew}

\begin{lem}
$S$ ist äquidimensional.
\end{lem}

\begin{bew}
Da $S$ ein \textsc{Cohen-Macaulay} lokaler Ring ist, folgt mit \cite[Thm. 17.4 (i)]{Mat}, dass für alle minimalen Primideale $\p$ von $S$ die Gleichung $\dim S / \p = \dim S$ erfüllt ist; somit ist $S$ äquidimensional nach Definition.
\end{bew}

\begin{lem}
$\dim S[1/2] = \dim S_\p = 8$ für alle maximalen Ideale $\p$ von $S[1/2]$.
\end{lem}

\begin{bew}
Da $S$ äquidimensional, als Quotient eines vollständigen lokalen noetherschen Rings vollständig, lokal und noethersch ist, und $2$ ein nicht nilpotentes Element im maximalen Ideal ist, folgt mit \cite[Lemma 2.3]{CDP} und Lemma \ref{vollstdurchschn}, dass $\dim S[1/2] = \dim S_\p = \dim S - 1 = 8$ für alle maximalen Ideale $\p$ von $S[1/2]$.
\end{bew}

\begin{prop}
\label{normal}
Die Ringe $S[1/2]$ und $S^\pm[1/2]$ sind normal. Insbesondere sind diese Ringe als (endliches) direktes Produkt von Integritätsbereichen darstellbar, und ihre irreduziblen Komponenten sind disjunkt im Spektrum.
\end{prop}

\begin{bew}
Nach der Bemerkung auf \cite[Seite 64]{Mat} und Lemma \ref{specspm12} genügt es zu zeigen, dass $S[1/2]$ normal ist.

Da $S[1/2]$ \textsc{Cohen-Macaulay} ist, erfüllt dieser Ring \textsc{Serre}'s Bedingung $S_2$. Wenn wir noch zeigen könnten, dass dieser Ring \textsc{Serre}'s Bedingung $R_1$ erfüllt, dann folgt mit \cite[Thm. 23.8]{Mat}, dass der noethersche Ring $S[1/2]$ normal ist.

Da $S[1/2]$ exzellent ist, folgt mit (3) in der Definition auf \cite[Seite 260]{Mat}, dass der singuläre Ort von $\Spec S[1/2]$ abgeschlossen ist; insbesondere gibt es ein Ideal $I$ von $S[1/2]$, so dass der singuläre Ort durch $V(I)$ gegeben ist. Nach Lemma \ref{singort} gilt $\dim S[1/2] / I \leq 6$.

Angenommen, es gibt ein $\p \in \Spec S[1/2]$ mit $\hgt \p \leq 1$ und $I \subseteq \p$. Dann sei $\q$ ein beliebiges maximales Ideal von $S[1/2]$, welches $\p$ enthält. Da $S_\q$ nach \cite[Seite 136]{Mat} ein \textsc{Cohen-Macaulay} lokaler Ring ist, folgt mit \cite[Thm. 17.4]{Mat}, dass $\dim S[1/2] / I \geq \dim S_\q / \p = \dim S_\q - \hgt \p \geq 8 - 1 = 7$, und wir erhalten einen Widerspruch. Somit gilt für alle $\p \in \Spec S[1/2]$ mit $\hgt \p \leq 1$, dass $I \nsubseteq \p$ bzw. $\p \notin V(I)$, d.~h. $S[1/2]_\p$ ist regulär, und $R_1$ ist erfüllt.

Zu "`Insbesondere"': Folgt mit \cite[Seite 64]{Mat}.
\end{bew}

\begin{lem}
\label{schnitt}
Sei $V \mathrel{\mathop:}= \Spec S[1/2] / (x_{11}, x_{12}, x_{21} + \tr \tilde X, y_{21}, z_{21})$. Dann wird jede irreduzible Komponente von $\Spec S[1/2]$ von $V$ geschnitten, und $V(K)$ ist die Menge der Tripel $(\tilde X, \tilde Y, \tilde Z)$ in
$$\left( \begin{pmatrix} 1 & \lambda \\ 0 & 1 \end{pmatrix} + \M_2(\m_K) \right) \times \left( \begin{pmatrix} 1 & \mu \\ 0 & 1 \end{pmatrix} + \M_2(\m_K) \right) \times \left( \begin{pmatrix} 1 & \kappa \\ 0 & 1 \end{pmatrix} + \M_2(\m_K) \right)$$
mit $\tilde X = \begin{pmatrix} 1 & \lambda \\ 0 & -1 \end{pmatrix}$, $\tilde Y = \begin{pmatrix} 1 + y_{11} & \mu + y_{12} \\ 0 & 1 + y_{22} \end{pmatrix}$, $\tilde Z = \begin{pmatrix} 1 + z_{11} & \kappa + z_{12} \\ 0 & 1 + z_{22} \end{pmatrix}$, so dass $(1 + y_{11})^4 = (1 + y_{22})^4 = 1$ und einer der folgenden drei Bedingungen erfüllt ist:
\begin{enumerate}[1.]
\item $y_{11} \neq y_{22}$ und $(y_{11} - y_{22}) (\kappa + z_{12}) = (\mu + y_{12}) (z_{11} - z_{22})$;
\item $y_{11} = y_{22}$ und $\mu = y_{12} = 0$; oder
\item $y_{11} = y_{22}$ und $(1 + z_{11}) = 5 (1 + z_{22})$.
\end{enumerate}
\end{lem}

\begin{bew}
Da $S$ lokal, noethersch und \textsc{Cohen-Macaulay} ist, und die Folge $\varpi, x_{11}, x_{12}, x_{21} + \tr \tilde X, y_{21}, z_{21}$ nach Lemma \ref{vollstdurchschn} regulär ist, folgt mit \cite[Prop. 5.1 b)]{CDP}, dass jede irreduzible Komponente von $\Spec S[1/2]$ von $V$ geschnitten wird. Die Elemente von $V(K)$ sind die Tripel $(\tilde X, \tilde Y, \tilde Z)$ in
$$\left( \begin{pmatrix} 1 & \lambda \\ 0 & 1 \end{pmatrix} + \M_2(\m_K) \right) \times \left( \begin{pmatrix} 1 & \mu \\ 0 & 1 \end{pmatrix} + \M_2(\m_K) \right) \times \left( \begin{pmatrix} 1 & \kappa \\ 0 & 1 \end{pmatrix} + \M_2(\m_K) \right),$$
so dass $\tilde X^2 \tilde Y^4 [\tilde Y, \tilde Z] = 1$ und $x_{11} = x_{12} = {x_{21} + \tr \tilde X} = y_{21} = z_{21} = 0$. Dies impliziert, dass $\tilde Y$, $\tilde Z$, und somit auch $\tilde Y^4 [\tilde Y, \tilde Z]$ und $\tilde X^2$ obere Dreiecksmatrizen sind; insbesondere haben wir $x_{21}^2 = -x_{21} \tr \tilde X = (\tr \tilde X)^2 = 0$, und somit ${x_{21} = \tr \tilde X = 0}$. Es folgt $\tilde X = \begin{pmatrix} 1 & \lambda \\ 0 & -1 \end{pmatrix}$ und $\tilde X^2 = 1$. Also vereinfacht sich die Relation $\tilde X^2 \tilde Y^4 [\tilde Y, \tilde Z] = 1$ zu $\tilde Y^4 [\tilde Y, \tilde Z] = 1$. Betrachten wir die Hauptdiagonalelemente, so erhalten wir $(1 + y_{11})^4 = (1 + y_{22})^4 = 1$. Nun unterscheiden wir zwei Fälle:

\textbf{(i) Es gilt $y_{11} \neq y_{22}$:} Teilen wir die Identität
$$(\y_{11} - \y_{22}) (\y_{11} + \y_{22}) (\y_{11}^2 + \y_{22}^2) = (\y_{11}^2 - \y_{22}^2) (\y_{11}^2 + \y_{22}^2) = \y_{11}^4 - \y_{22}^4 = 1 - 1 = 0$$
durch $\y_{11} - \y_{22} = y_{11} - y_{22} \neq 0$, so erhalten wir $(\y_{11} + \y_{22}) (\y_{11}^2 + \y_{22}^2) = 0$, und somit $\tilde Y^4 = \begin{pmatrix} \y_{11}^4 & \y_{12} (\y_{11} + \y_{22}) (\y_{11}^2 + \y_{22}^2) \\ 0 & \y_{22}^4 \end{pmatrix} = 1$; also vereinfacht sich die Relation zu $[\tilde Y, \tilde Z] = 1$ bzw. $\tilde Y \tilde Z = \tilde Z \tilde Y$. Dies impliziert, dass die Determinanten aller $2 \times 2$-Matrizen in $\begin{pmatrix} y_{11} - y_{22} & \mu + y_{12} & 0 \\ z_{11} - z_{22} & \kappa + z_{12} & 0 \end{pmatrix}$ verschwinden, d.~h. wir haben $(y_{11} - y_{22}) (\kappa + z_{12}) = (\mu + y_{12}) (z_{11} - z_{22})$, und landen in 1.

\textbf{(ii) Es gilt $y_{11} = y_{22}$:} Dann erhalten wir durch Ausrechnen, dass
\begin{align*}
\tilde Y^5 \tilde Z &= \begin{pmatrix} 1 + y_{11} & 5 (\mu + y_{12}) \\ 0 & 1 + y_{11} \end{pmatrix} \begin{pmatrix} 1 + z_{11} & \kappa + z_{12} \\ 0 & 1 + z_{22} \end{pmatrix} \\ &= \begin{pmatrix} (1 + y_{11}) (1 + z_{11}) & (1 + y_{11}) (\kappa + z_{12}) + 5 (\mu + y_{12}) (1 + z_{22}) \\ 0 & (1 + y_{11}) (1 + z_{22}) \end{pmatrix}
\end{align*}
und
$$\tilde Z \tilde Y = \begin{pmatrix} (1 + y_{11}) (1 + z_{11}) & (\mu + y_{12}) (1 + z_{11}) + (1 + y_{11}) (\kappa + z_{12}) \\ 0 & (1 + y_{11}) (1 + z_{22}) \end{pmatrix}$$
gelten. Da die Relation  $\tilde Y^4 [\tilde Y, \tilde Z] = 1$ zu $\tilde Y^5 \tilde Z = \tilde Z \tilde Y$ äquivalent ist, erhalten wir
$$(1 + y_{11}) (\kappa + z_{12}) + 5 (\mu + y_{12}) (1 + z_{22}) = (\mu + y_{12}) (1 + z_{11}) + (1 + y_{11}) (\kappa + z_{12})$$
bzw. $(\mu + y_{12}) (5 (1 + z_{22}) - (1 + z_{11})) = 0,$ was genau dann der Fall ist, wenn $\mu + y_{12} = 0$, also $\mu = y_{12} = 0$, oder aber $(1 + z_{11}) = 5 (1 + z_{22})$ erfüllt ist, und landen in 2. oder 3.
\end{bew}

\begin{bem}
\label{regfolge}
Die reguläre Folge $\varpi, x_{11}, x_{12}, x_{21} + \tr \tilde X, y_{21}, z_{21}$ ist anders gewählt als in \cite{CDP} und stellt eine technische Verbesserung der Arbeit \cite{CDP} dar.
\end{bem}

\newpage

\section{Punkte auf 
Komponenten von \texorpdfstring{$S[1/2]$}{S[1/2]}}

\setcounter{prop}{-1}

Wir werden zeigen, dass es Punkte bestimmter Form auf den irreduziblen Komponenten von $S[1/2]$ gibt.

\begin{bem}
Ein regulärer lokaler Ring ist nach \cite[Thm. 20.3]{Mat} stets faktoriell.
\end{bem}

\begin{lem}
\label{intber1}
Sei $A$ ein regulärer lokaler Ring mit maximalem Ideal $\m$; insbesondere ist dann auch $A[[x]]$ nach \cite[Thm. 20.8]{Mat} ein regulärer faktorieller Ring. Ist $a \in A$, $b \in \m$ und entweder $a \in A^\times$ oder $b$ prim in $A$, so ist $A[[x]] / (ax + b)$ ein Integritätsbereich.
\end{lem}

\begin{bew}
Angenommen, es gelte $ax + b = (a_0 + a_1 x + \ldots) (b_0 + b_1 x + \ldots)$, wobei $a_0, b_0 \in \m$. Dann ist klar, das wir $a_0 b_0 = b$ und $a_0 b_1 + a_1 b_0 = a$ haben; insbesondere folgt daraus, dass $b$ kein Primelement in $A$ sein kann, und dass wir $a \in \m$ haben, was aber im Widerspruch zu unseren Voraussetzungen steht. Also ist $ax + b$ in $A[[x]]$ irreduzibel und damit prim, da $A[[x]]$ ein faktorieller Ring ist, und daraus folgt die Behauptung.
\end{bew}

\begin{lem}
\label{intber2}
Sei $A$ ein regulärer lokaler Ring mit maximalem Ideal $\m$; insbesondere ist dann auch $A[[x, y]]$ nach \cite[Thm. 20.8]{Mat} ein regulärer faktorieller Ring. Ist $a \in A^\times$ und $b \in \m \setminus \{0\}$, so ist $A[[x, y]] / (axy + b)$ ein Integritätsbereich.
\end{lem}

\begin{bew}
Angenommen, es gelte 
$$axy + b = (a_{00} + a_{10} x + a_{01} y + a_{11} xy + \ldots) (b_{00} + b_{10} x + b_{01} y + b_{11} xy + \ldots),$$
wobei $a_{00}, b_{00} \in \m$. Dann ist klar, dass wir
$$a_{00} b_{00} = b, \quad a_{00} b_{10} + a_{10} b_{00} = 0,\quad 
a_{00} b_{01} + a_{01} b_{00} = 0,$$
$$a_{00} b_{20} + a_{10} b_{10} + a_{20} b_{00} = 0, \quad 
a_{00} b_{11} + a_{10} b_{01} + a_{01} b_{10} + a_{11} b_{00} = a$$
haben. Daraus folgen $a_{10} b_{10} = -(a_{00} b_{20} + a_{20} b_{00}) \in \m$ und
\begin{align*}
a_{00} (a_{10} b_{01} - a_{01} b_{10}) &= a_{00} a_{10} b_{01} - a_{00} a_{01} b_{10} \\ &= a_{10} a_{00} b_{01} + a_{10} a_{01} b_{00} - a_{01} a_{00} b_{10} - a_{01} a_{10} b_{00} \\ &= a_{10} (a_{00} b_{01} + a_{01} b_{00}) - a_{01} (a_{00} b_{10} + a_{10} b_{00}) = 0;
\end{align*}
somit folgen $a_{10} \in \m$ oder $b_{10} \in \m$ und nach Division durch $a_{00} \neq 0$ die Gleichung $a_{10} b_{01} = a_{01} b_{10} \in \m$. Es folgt $a = a_{00} b_{11} + a_{10} b_{01} + a_{01} b_{10} + a_{11} b_{00} \in \m$, im Widerspruch zu $a \in A^\times$. Also ist $axy + b$ in $A[[x, y]]$ irreduzibel und damit prim, da $A[[x, y]]$ ein faktorieller Ring ist, und daraus folgt die Behauptung. 
\end{bew}

\begin{satz}
\label{4irredkomp}
Der Ring $S[1/2] / (x_{21}, y_{21}, z_{21})$ hat $4$ irreduzible Komponenten, die durch die Gleichungen $\x_{11} \y_{11}^2 = \varepsilon_1$, $\x_{22} \y_{22}^2 = \varepsilon_2$ für $\varepsilon_1, \varepsilon_2 \in \{\pm1\}$ gegeben sind.
\end{satz}

\begin{bew}
Durch Ausmultiplizieren erhalten wir
$$\tilde X^2 \tilde Y^4 = \begin{pmatrix} \x_{11}^2 \y_{11}^4 & \x_{11}^2 \y_{12} (\y_{11} + \y_{22}) (\y_{11}^2 + \y_{22}^2) + \x_{12} (\x_{11} + \x_{22}) \y_{22}^4 \\ 0 & \x_{22}^2 \y_{22}^4 \end{pmatrix},$$
$$[\tilde Y, \tilde Z] =  \begin{pmatrix} 1 & ((\y_{11} - \y_{22}) \z_{12} + \y_{12} (\z_{22} - \z_{11})) \y_{22}^{-1} \z_{22}^{-1} \\ 0 & 1 \end{pmatrix}.$$
Damit ist die Relation $\tilde X^2 \tilde Y^4 [\tilde Y, \tilde Z] = 1$ äquivalent zu $\x_{11}^2 \y_{11}^4 = 1$, $\x_{22}^2 \y_{22}^4 = 1$ und
\begin{displaymath}
\begin{split}
\x_{11}^2 \y_{12} (\y_{11} + \y_{22}) (\y_{11}^2 + \y_{22}^2) &+ \x_{12} (\x_{11} + \x_{22}) \y_{22}^4 \\ &+ \frac{(\y_{11} - \y_{22}) \z_{12} + \y_{12} (\z_{22} - \z_{11})}{\y_{22} \z_{22}} = 0.
\end{split}
\end{displaymath}
Die Relationen $\x_{11}^2 \y_{11}^4 = 1$ und $\x_{22}^2 \y_{22}^4 = 1$ implizieren, dass wenn $\p$ eine irreduzible Komponente von $\Spec S[1/2] / (x_{21}, y_{21}, z_{21})$ ist, dann sind $\x_{11} \y_{11}^2 = \varepsilon_1$ und $\x_{22} \y_{22}^2 = \varepsilon_2$ in $S[1/2] / (x_{21}, y_{21}, z_{21}, \p)$ für gewisse $\varepsilon_1, \varepsilon_2 \in \{\pm1\}$ erfüllt.

Also reicht es zu zeigen, dass $S[1/2] / (x_{21}, y_{21}, z_{21}, \x_{11} \y_{11}^2 - \varepsilon_1, \x_{22} \y_{22}^2 - \varepsilon_2)$ für $\varepsilon_1, \varepsilon_2 \in \{\pm1\}$ ein Integritätsbereich ist.

Die Relationen $\x_{11} \y_{11}^2 = \varepsilon_1$ und $\x_{22} \y_{22}^2 = \varepsilon_2$ implizieren, dass $\x_{11} = \varepsilon_1 \y_{11}^{-2}$ und $\x_{22} = \varepsilon_2 \y_{22}^{-2}$. Somit vereinfacht sich die Relation $\tilde X^2 \tilde Y^4 [\tilde Y, \tilde Z] = 1$ zu 
\begin{displaymath}
\begin{split}
\y_{11}^{-1} \y_{12} (1 + \y_{11}^{-1} \y_{22}) (1 + \y_{11}^{-2} \y_{22}^2) &+ \x_{12} (\varepsilon_1 \y_{11}^{-2} \y_{22}^2 +\varepsilon_2) \y_{22}^2 \\ &+ \frac{(\y_{11} - \y_{22}) \z_{12} + \y_{12} (\z_{22} - \z_{11})}{\y_{22} \z_{22}} = 0.
\end{split}
\end{displaymath}
Erweitern mit $\y_{11}^4 \y_{22} \z_{22}$ liefert 
\begin{displaymath}
\begin{split}
\y_{12} (\y_{11} + \y_{22}) (\y_{11}^2 + \y_{22}^2) \y_{22} \z_{22} &+ \x_{12} \y_{11}^2 (\varepsilon_1 \y_{22}^2 +\varepsilon_2 \y_{11}^2) \y_{22}^3 \z_{22} \\ &+ ((\y_{11} - \y_{22}) \z_{12} + \y_{12} (\z_{22} - \z_{11})) \y_{11}^4 = 0.
\end{split}
\end{displaymath}
Nach Definition von $\tilde X$, $\tilde Y$ und $\tilde Z$ haben wir $\y_{11} = 1 + y_{11}$, $\y_{22} = 1 + y_{22}$, $\z_{11} = 1 + z_{11}$, $\z_{22} = 1 + z_{22}$, $\x_{12} = \lambda + x_{12}$, $\y_{12} = \mu + y_{12}$ und $\z_{12} = \kappa + z_{12}$, und müssen zeigen, dass 
\begin{displaymath}
\begin{split}
f \mathrel{\mathop:}= \y_{12} (\y_{11} + \y_{22}) (\y_{11}^2 + \y_{22}^2) \y_{22} \z_{22} &+ \x_{12} \y_{11}^2 (\varepsilon_1 \y_{22}^2 + \varepsilon_2 \y_{11}^2) \y_{22}^3 \z_{22} \\ &+ ((\y_{11} - \y_{22}) \z_{12} + \y_{12} (\z_{22} - \z_{11})) \y_{11}^4
\end{split}
\end{displaymath}
ein Primelement in $\Oc[[x_{12}, y_{12}, z_{12}, y_{11}, y_{22}, z_{11}, z_{22}]][1/2]$ ist.

Angenommen, das Element $f$ ist in $\Oc[[x_{12}, y_{12}, z_{12}, y_{11}, y_{22}, z_{11}, z_{22}]]$ reduzibel. Wir unterscheiden zwei Fälle:

\textbf{(i) Es gilt $\mu \neq 0$:} Wenn wir $y_{22}$ für $y_{11}$, also $\y_{22}$ für $\y_{11}$ in $f$ einsetzen, so erhalten wir
$$(\varepsilon_1 + \varepsilon_2) \x_{12} \y_{22}^7 \z_{22} + \y_{12} (5 \z_{22} - \z_{11}) \y_{22}^4,$$
und dieses Element ist dann in $\Oc[[x_{12}, y_{12}, z_{12}, y_{22}, z_{11}, z_{22}]]$ reduzibel. Division durch die Einheit $\y_{22}^4$ liefert
$$(\varepsilon_1 + \varepsilon_2) \x_{12} \y_{22}^3 \z_{22} + \y_{12} (5 \z_{22} - \z_{11}),$$
und dieses Element ist ebenfalls in $\Oc[[x_{12}, y_{12}, z_{12}, y_{22}, z_{11}, z_{22}]]$ reduzibel. Durch die Relation $5 \z_{22} - \z_{11} = z$ erhalten wir einen Isomorphismus
$$\Oc[[x_{12}, y_{12}, z_{12}, y_{22}, z_{11}, z_{22}]] \cong \Oc[[x_{12}, y_{12}, z_{12}, y_{22}, z_{22}, z]];$$
damit wäre das Element $(\varepsilon_1 + \varepsilon_2) \x_{12} \y_{22}^3 \z_{22} + \y_{12} z$ in $\Oc[[x_{12}, y_{12}, z_{12}, y_{22}, z_{22}, z]]$ reduzibel.

Nun folgt jedoch mit der Tatsache, dass $(\varepsilon_1 + \varepsilon_2) \x_{12} \y_{22}^3 \z_{22}$ im maximalen Ideal von $\Oc[[x_{12}, y_{12}, z_{12}, y_{22}, z_{22}]]$ liegt, und dass wegen $\mu \neq 0$ das Element $\y_{12}$ eine Einheit in $\Oc[[x_{12}, y_{12}, z_{12}, y_{22}, z_{22}]]$ ist, dass das Element $(\varepsilon_1 + \varepsilon_2) \x_{12} \y_{22}^3 \z_{22} + \y_{12} z$ in $\Oc[[x_{12}, y_{12}, z_{12}, y_{22}, z_{22}, z]]$ nach Lemma \ref{intber1} irreduzibel ist. Wir erhalten somit einen Widerspruch.

\textbf{(ii) Es gilt $\mu = 0$:} Falls $\varepsilon_1 = \varepsilon_2 \in \{ 1, -1 \}$ gilt, erhalten wir wie im Fall (i), dass das Element $(\varepsilon_1 + \varepsilon_2) \x_{12} \y_{22}^3 \z_{22} + \y_{12} z$ in $\Oc[[x_{12}, y_{12}, z_{12}, y_{22}, z_{22}, z]]$ reduzibel ist. Da $\mu = 0$, kann dieses Element auch als $(\varepsilon_1 + \varepsilon_2) \x_{12} \y_{22}^3 \z_{22} + y_{12} z$ geschrieben werden.

Nun folgt jedoch mit der Tatsache, dass $(\varepsilon_1 + \varepsilon_2) \x_{12} \y_{22}^3 \z_{22}$ ein von 0 verschiedenes Element im maximalen Ideal von $\Oc[[x_{12}, z_{12}, y_{22}, z_{22}]]$ ist, dass das Element $(\varepsilon_1 + \varepsilon_2) \x_{12} \y_{22}^3 \z_{22} + y_{12} z$ in $\Oc[[x_{12}, y_{12}, z_{12}, y_{22}, z_{22}, z]]$ nach Lemma \ref{intber2} irreduzibel ist. Wir erhalten somit einen Widerspruch.

Ansonsten gilt $\varepsilon_1 + \varepsilon_2 = 0$, und mit Einsetzen von $-y_{22} - 2$ für $y_{11}$, also $-\y_{22}$ für $\y_{11}$ erhalten wir $(-2 \y_{22} \z_{12} + \y_{12} (\z_{22} - \z_{11})) \y_{22}^4$, und dieses Element ist dann in $\Oc[[x_{12}, y_{12}, z_{12}, y_{22}, z_{11}, z_{22}]]$ reduzibel. Division durch die Einheit $\y_{22}^4$ liefert $-2 \y_{22} \z_{12} + \y_{12} (\z_{22} - \z_{11})$, und dieses Element ist ebenfalls in $\Oc[[x_{12}, y_{12}, z_{12}, y_{22}, z_{11}, z_{22}]]$ reduzibel. Durch die Relation $\z_{22} - \z_{11} = z$ erhalten wir einen Isomorphismus
$$\Oc[[x_{12}, y_{12}, z_{12}, y_{22}, z_{11}, z_{22}]] \cong \Oc[[x_{12}, y_{12}, z_{12}, y_{22}, z_{22}, z]];$$
damit wäre das Element $-2 \y_{22} \z_{12} + \y_{12} z$ in $\Oc[[x_{12}, y_{12}, z_{12}, y_{22}, z_{22}, z]]$ reduzibel. Da $\mu = 0$, kann dieses Element auch als $-2 \y_{22} \z_{12} + y_{12} z$ geschrieben werden.

Nun folgt jedoch mit der Tatsache, dass $-2 \y_{22} \z_{12}$ ein von $0$ verschiedenes Element im maximalen Ideal von $\Oc[[x_{12}, z_{12}, y_{22}, z_{22}]]$ ist, dass das Element $(\varepsilon_1 + \varepsilon_2) \x_{12} \y_{22}^3 \z_{22} + y_{12} z$ in $\Oc[[x_{12}, y_{12}, z_{12}, y_{22}, z_{22}, z]]$ nach Lemma \ref{intber2} irreduzibel ist. Wir erhalten somit einen Widerspruch.

Mit (i) und (ii) folgt nun, dass das genannte Element $f$ im faktoriellen Ring $\Oc[[x_{12}, y_{12}, z_{12}, y_{11}, y_{22}, z_{11}, z_{22}]]$ irreduzibel und somit prim ist. Da auch die Lokalisierung $\Oc[[x_{12}, y_{12}, z_{12}, y_{11}, y_{22}, z_{11}, z_{22}]][1/2]$ ein faktorieller Ring ist, und $\varpi$ kein Teiler von $f$ ist (der Koeffizient von $y_{12} z_{11}$ in $f$ ist gleich $-1$), ist $f$ auch in $\Oc[[x_{12}, y_{12}, z_{12}, y_{11}, y_{22}, z_{11}, z_{22}]][1/2]$ ein Primelement, und die Behauptung folgt.
\end{bew}

\begin{kor}
\label{punkte1}
Sei $K$ eine endliche Körpererweiterung von $\Q_2$, welche ein Element $\zeta_8$ mit $\zeta_8^4 = -1$ enthält; ferner sei $i = \zeta_8^2$. Falls $\mu = 0$, so enthält jede irreduzible Komponente von $\Spec S[1/2]$ einen der $K$-Punkte 

$\left( \tilde X_\lambda \mathrel{\mathop:}= \begin{pmatrix} 1 & \lambda \\ 0 & -1 \end{pmatrix}, \tilde Y_{\kappa, n} \mathrel{\mathop:}= \tilde Z_{\kappa, n}^2, \tilde Z_{\kappa, n} \right) \quad (1 \leq n \leq 4)$,

wobei

$\tilde Z_{\kappa, 1} \mathrel{\mathop:}= \begin{pmatrix} 1 & \kappa \\ 0 & \zeta_8 \end{pmatrix}$, $\tilde Z_{\kappa, 2} \mathrel{\mathop:}= \begin{pmatrix} 1 & \kappa \\ 0 & i \end{pmatrix}$, $\tilde Z_{\kappa, 3} \mathrel{\mathop:}= \begin{pmatrix} \zeta_8 & \kappa \\ 0 & i \end{pmatrix}$, und $\tilde Z_{\kappa, 4} \mathrel{\mathop:}= \begin{pmatrix} \zeta_8 & \kappa \\ 0 & \zeta_8^3 \end{pmatrix}$.
\end{kor}

\begin{bew}
Nach Proposition \ref{normal} wissen wir, dass die irreduziblen Komponenten von $\Spec S[1/2]$ disjunkt sind. Nach Lemma \ref{schnitt} enthält jede irreduzible Komponente von $\Spec S[1/2]$ Punkte mit $x_{21} = y_{21} = z_{21} = 0$. Außerdem folgt mit Satz \ref{4irredkomp}, dass $\Spec S[1/2] / (x_{21}, y_{21}, z_{21})$ vier irreduzible Komponenten hat; jede dieser irreduziblen Komponenten ist in einer solchen von $\Spec S[1/2]$ enthalten.

Man kann nachprüfen, dass $\tilde X_\lambda$, $\tilde Y_{\kappa, n}$ und $\tilde Z_{\kappa, n}$ modulo $\m_K$ sich zu $\begin{pmatrix} 1 & \bar\lambda \\ 0 & 1 \end{pmatrix}$, $\begin{pmatrix} 1 & 0 \\ 0 & 1 \end{pmatrix}$ und $\begin{pmatrix} 1 & \bar\kappa \\ 0 & 1 \end{pmatrix}$ reduzieren.

Ferner haben wir $\tilde X_\lambda^2 = 1$, $\tilde Y_{\kappa, n}^4 = \tilde Z_{\kappa, n}^8 = 1$ und $[\tilde Y_{\kappa, n}, \tilde Z_{\kappa, n}] = [\tilde Z_{\kappa, n}^2, \tilde Z_{\kappa, n}] = 1$, woraus $\tilde X_\lambda^2 \tilde Y_{\kappa, n}^4 [\tilde Y_{\kappa, n}, \tilde Z_{\kappa, n}] = 1$ für $1 \leq n \leq 4$ folgt.

Da $\tilde X_\lambda$ und die $\tilde Y_{\kappa, n}, \tilde Z_{\kappa, n}$ obere Dreiecksmatrizen sind, erhalten wir durch Ausrechnen der Elemente $\varepsilon_1 = \x_{11} \y_{11}^2 = \z_{11}^4, \varepsilon_2 = \x_{22} \y_{22}^2 = -\z_{22}^4$, dass jeder der vier genannten $K$-Punkte zu einer anderen irreduziblen Komponente von $\Spec S[1/2] / (x_{21}, y_{21}, z_{21})$ gehört. Dies impliziert, dass jede irreduzible Komponente von $\Spec S[1/2]$ mindestens einen dieser vier $K$-Punkte enthält.
\end{bew}

\begin{kor}
\label{punkte2}
Sei $K$ eine endliche Körpererweiterung von $\Q_2$, welche ein Element $i$ mit $i^2 = -1$ enthält. Falls $\mu \neq 0$, so enthält jede irreduzible Komponente von $\Spec S[1/2]$ einen der $K$-Punkte 

$\left( \tilde X_\lambda \mathrel{\mathop:}= \begin{pmatrix} 1 & \lambda \\ 0 & -1 \end{pmatrix}, \tilde Y_{\mu, n}, \tilde Z_{\mu\kappa, n} \mathrel{\mathop:}= 1 + \frac\kappa\mu (\tilde Y_{\mu, n} - 1) \right) \quad (1 \leq n \leq 4)$,

wobei

$\tilde Y_{\mu, 1} \mathrel{\mathop:}= \begin{pmatrix} 1 & \mu \\ 0 & i \end{pmatrix}$, $\tilde Y_{\mu, 2} \mathrel{\mathop:}= \begin{pmatrix} 1 & \mu \\ 0 & -1 \end{pmatrix}$, $\tilde Y_{\mu, 3} \mathrel{\mathop:}= \begin{pmatrix} i & \mu \\ 0 & -1 \end{pmatrix}$, und $\tilde Y_{\mu, 4} \mathrel{\mathop:}= \begin{pmatrix} i & \mu \\ 0 & -i \end{pmatrix}$.
\end{kor}

\begin{bew}
Nach Proposition \ref{normal} wissen wir, dass die irreduziblen Komponenten von $\Spec S[1/2]$ disjunkt sind. Nach Lemma \ref{schnitt} enthält jede irreduzible Komponente von $\Spec S[1/2]$ Punkte mit $x_{21} = y_{21} = z_{21} = 0$. Außerdem folgt mit Satz \ref{4irredkomp}, dass $\Spec S[1/2] / (x_{21}, y_{21}, z_{21})$ vier irreduzible Komponenten hat; jede dieser irreduziblen Komponenten ist in einer solchen von $\Spec S[1/2]$ enthalten.

Man kann nachprüfen, dass $\tilde X_\lambda$, $\tilde Y_{\mu, n}$ und $\tilde Z_{\mu\kappa, n}$ modulo $\m_K$ sich zu $\begin{pmatrix} 1 & \bar\lambda \\ 0 & 1 \end{pmatrix}$, $\begin{pmatrix} 1 & \bar\mu \\ 0 & 1 \end{pmatrix}$ und $\begin{pmatrix} 1 & \bar\kappa \\ 0 & 1 \end{pmatrix}$ reduzieren.

Ferner haben wir $\tilde X_\lambda^2 = 1$, $\tilde Y_{\mu, n}^4 = 1$ und $$[\tilde Y_{\mu, n}, \tilde Z_{\mu\kappa, n}] = [\tilde Y_{\mu, n}, 1 + \frac\kappa\mu (\tilde Y_{\mu, n} - 1)] = 1,$$ woraus $\tilde X_\lambda^2 \tilde Y_{\mu, n}^4 [\tilde Y_{\mu, n}, \tilde Z_{\mu\kappa, n}] = 1$ für $1 \leq n \leq 4$ folgt.

Da $\tilde X_\lambda$ und die $\tilde Y_{\mu, n}, \tilde Z_{\mu\kappa, n}$ obere Dreiecksmatrizen sind, erhalten wir durch Ausrechnen der Elemente $\varepsilon_1 = \x_{11} \y_{11}^2 = \y_{11}^2, \varepsilon_2 = \x_{22} \y_{22}^2 = -\y_{22}^2$, dass jeder der vier genannten $K$-Punkte zu einer anderen irreduziblen Komponente von $\Spec S[1/2] / (x_{21}, y_{21}, z_{21})$ gehört. Dies impliziert, dass jede irreduzible Komponente von $\Spec S[1/2]$ mindestens einen dieser vier $K$-Punkte enthält.
\end{bew}


\section{\texorpdfstring{$S^+[1/2]$}{S+[1/2]} und \texorpdfstring{$S^-[1/2]$}{S-[1/2]} sind Integritätsbereiche}
\label{disks}

In diesem Kapitel zeigen wir mit der Methode von \cite{CDP}, dass $S^+[1/2]$ und $S^-[1/2]$ Integritätsbereiche sind.

\begin{defi}
Sei $A$ ein vollständige lokale noethersche $\Oc$-Algebra mit Restklassenkörper $k$, und sei $X = \Spec A[1/2]$. Außerdem sei $K$ eine beliebige endliche Körpererweiterung von $L$, $\Oc_K$ der Ganzheitsring von $K$ und $\m_K$ das maximale Ideal von $\Oc_K$. Sei $T_K$ die \textsc{Tate}-Algebra in einer Variablen über $K$, d.~h. der Ring der Potenzreihen in $\Oc_K[[t]][1/2]$, welche auf ganz $\Oc_{\C_2}$ konvergieren, wobei $\C_2$ die $2$-adische Vervollständigung des algebraischen Abschlusses von $\Q_2$ ist.
\end{defi}

\begin{defi}
\label{bogenkette}
Wir sagen, dass $x_0, x_1 \in X(K)$ durch einen Bogen verbunden sind, falls es einen $\Oc$-Algebrenhomomorphismus $\varphi : A \rightarrow T_K$ gibt, so dass $x_0 : A \rightarrow K$ durch Spezialisierung von $\varphi$ bei $t = 0$ und $x_1 : A \rightarrow K$ durch Spezialisierung von $\varphi$ bei $t = 1$ erhalten wird. Wir sagen, dass $x_0, x_1 \in X(K)$ durch eine Kette verbunden sind, falls es eine endliche Folge von Elementen von $X(K)$, nämlich $x_0 = y_0, y_1, \ldots , y_n = x_1$, gibt, so dass $y_{i-1}, y_i$ für alle $1 \leq i \leq n$ durch einen Bogen verbunden sind.
\end{defi}

\begin{beisp}
Für $A = S$ steht die Menge $X(K)$ in Bijektion zur Menge der Tripel $(\tilde X, \tilde Y, \tilde Z)$ der Matrizen in
$$\left( \begin{pmatrix} 1 & \lambda \\ 0 & 1 \end{pmatrix} + \M_2(\m_K) \right) \times \left( \begin{pmatrix} 1 & \mu \\ 0 & 1 \end{pmatrix} + \M_2(\m_K) \right) \times \left( \begin{pmatrix} 1 & \kappa \\ 0 & 1 \end{pmatrix} + \M_2(\m_K) \right)$$
mit $\tilde X^2 \tilde Y^4 [\tilde Y, \tilde Z] = 1$. Um zu zeigen, dass zwei Punkte, die zu den Tripeln $(\tilde X_0, \tilde Y_0, \tilde Z_0)$ und $(\tilde X_1, \tilde Y_1, \tilde Z_1)$ korrespondieren, durch einen Bogen verbunden sind, reicht es, Matrizen $(\tilde X(t), \tilde Y(t), \tilde Z(t))$ in $\M_2(T_K)$ zu produzieren, so dass folgende Bedingungen erfüllt sind:

\begin{itemize}
\item Die Einträge von $\tilde X - \begin{pmatrix} 1 & \lambda \\ 0 & 1 \end{pmatrix}$, $\tilde Y -  \begin{pmatrix} 1 & \mu \\ 0 & 1 \end{pmatrix}$ und $\tilde Z - \begin{pmatrix} 1 & \kappa \\ 0 & 1 \end{pmatrix}$ sind topologisch nilpotent in $T_K$;
\item $\tilde X(t)^2 \tilde Y(t)^4 [\tilde Y(t), \tilde Z(t)] = 1$;
\item $(\tilde X_0, \tilde Y_0, \tilde Z_0) = (\tilde X(0), \tilde Y(0), \tilde Z(0))$, $(\tilde X_1, \tilde Y_1, \tilde Z_1) = (\tilde X(1), \tilde Y(1), \tilde Z(1))$.
\end{itemize}
\end{beisp}

\begin{lem}
\label{bogen}
Falls $x, y \in X(K)$ durch einen Bogen verbunden sind, dann liegen sie auf der gleichen irreduziblen Komponente von $X$.
\end{lem}

\begin{bew}
Sei $\varphi: A \rightarrow T_K$ ein Bogen, der $x$ und $y$ verbindet. Der Kern von $\varphi$ ist ein Primideal von $A$, weil $T_K$ ein Integritätsbereich ist. Damit enthält $\ker \varphi$ ein minimales Primideal $\q$ von $A[1/2]$, und sowohl $x$ als auch $y$ liegt auf $V(\q)$.
\end{bew}

\begin{lem}
Wir nehmen an, dass die irreduziblen Komponenten von $X$ disjunkt im Spektrum sind. Falls $x, y \in X(K)$ durch eine Kette verbunden sind, dann liegen sie auf der gleichen irreduziblen Komponente von $X$.
\end{lem}

\begin{bew}
Die Annahme impliziert, dass jedes $x \in X(K)$ auf einer eindeutigen irreduziblen Komponente von $X$ liegt. Die Behauptung folgt mit Lemma \ref{bogen}.
\end{bew}

\begin{lem}
\label{bogen1}
Falls $\mu = 0$, so sind unter den Voraussetzungen von Korollar \ref{punkte1} die $K$-Punkte $\left( \tilde X_\lambda, \tilde Y_{\kappa, 1}, \tilde Z_{\kappa, 1} \right)$ und $\left( \tilde X_\lambda, \tilde Y_{\kappa, 3}, \tilde Z_{\kappa, 3} \right)$ in $\Spec S^-[1/2]$, und $\left( \tilde X_\lambda, \tilde Y_{\kappa, 2}, \tilde Z_{\kappa, 2} \right)$ und $\left( \tilde X_\lambda, \tilde Y_{\kappa, 4}, \tilde Z_{\kappa, 4} \right)$ in $\Spec S^+[1/2]$ jeweils durch einen Bogen verbunden.
\end{lem}

\begin{bew}
Sei $a(t) \mathrel{\mathop:}= 1 + (\zeta_8 - 1) t$. Dann benutzen wir den Bogen
$$\left( \tilde X(t), \quad \tilde Y_n(t) = \tilde Z_n(t)^2, \quad \tilde Z_n(t) \right),$$
wobei
$$\tilde X(t) \mathrel{\mathop:}= \begin{cases} a(t)^{-4} \begin{pmatrix} (1 - 6t^2 + 4t^3) & t (1 - t) (2 + 4t) \\ t (1 - t) (6 - 4t) & -(1 - 6t^2 + 4t^3) \end{pmatrix}, &\text{ falls } \lambda = 0 \\ \begin{pmatrix} (1 - 2t) a(t)^{-4} & \lambda \\ \frac4\lambda t (1 - t) a(t)^{-8} & -(1 - 2t) a(t)^{-4} \end{pmatrix}, &\text{ falls } \lambda \neq 0 \end{cases},$$
$\tilde Z_1(t) \mathrel{\mathop:}= \begin{pmatrix} a(t) & \kappa \\ 0 & \zeta_8 a(t) \end{pmatrix}$ und $\tilde Z_2(t) \mathrel{\mathop:}= \begin{pmatrix} a(t) & \kappa \\ 0 & i a(t) \end{pmatrix}$.

Man kann nachprüfen, dass $\tilde X(t)$, $\tilde Y_n(t)$ und $\tilde Z_n(t)$ modulo $\m_K$ sich zu $\begin{pmatrix} 1 & \bar\lambda \\ 0 & 1 \end{pmatrix}$, $\begin{pmatrix} 1 & 0 \\ 0 & 1 \end{pmatrix}$ und $\begin{pmatrix} 1 & \bar\kappa \\ 0 & 1 \end{pmatrix}$ reduzieren. Durch Ausrechnen erhalten wir
$$\tilde X(t)^2 = a(t)^{-8}, \tilde Y_n(t)^4 = \tilde Z_n(t)^8 = a(t)^8$$
und
$$[\tilde Y_n(t), \tilde Z_n(t)] = [\tilde Z_n(t)^2, \tilde Z_n(t)] = 1,$$
woraus
$$\tilde X(t)^2 \tilde Y_n(t)^4 [\tilde Y_n(t), \tilde Z_n(t)] = 1 \text{ für } 1 \leq n \leq 2$$
folgt. Ferner haben wir $\delta_n = \det \tilde X(t) (\det \tilde Y_n(t))^2 = -a(t)^{-8} (\det \tilde Z_n(t))^4$, also $\delta_1 = -a(t)^{-8} (\zeta_8 a(t)^2)^4 = 1$ und $\delta_2 = -a(t)^{-8} (i a(t)^2)^4 = -1$. Somit gehört der Bogen für $n = 1$ zu $\Spec S^-[1/2]$ und für $n = 2$ zu $\Spec S^+[1/2]$.

Außerdem haben wir $a(0) = 1$ und $a(1) = \zeta_8$, somit folgen $\tilde X(0) = \tilde X(1) = \tilde X_\lambda$, $\tilde Z_n(0) = \tilde Z_{\kappa, n}$ und $\tilde Z_n(1) = \tilde Z_{\kappa, n + 2}$, ferner $\tilde Y_n(0) = \tilde Z_n(0)^2 = \tilde Z_{\kappa, n}^2 = \tilde Y_{\kappa, n}$ und $\tilde Y_n(1) = \tilde Z_n(1)^2 = \tilde Z_{\kappa, n + 2}^2 = \tilde Y_{\kappa, n + 2}$. Damit folgt die Behauptung.
\end{bew}

\begin{lem}
\label{bogen2}
Falls $\mu \neq 0$, so sind unter den Voraussetzungen von Korollar \ref{punkte2} die $K$-Punkte $\left( \tilde X_\lambda, \tilde Y_{\mu, 1}, \tilde Z_{\mu\kappa, 1} \right)$ und $\left( \tilde X_\lambda, \tilde Y_{\mu, 3}, \tilde Z_{\mu\kappa, 3} \right)$ in $\Spec S^-[1/2]$, und $\left( \tilde X_\lambda, \tilde Y_{\mu, 2}, \tilde Z_{\mu\kappa, 2} \right)$ und $\left( \tilde X_\lambda, \tilde Y_{\mu, 4}, \tilde Z_{\mu\kappa, 4} \right)$ in $\Spec S^+[1/2]$ jeweils durch einen Bogen verbunden (für eine Definition des Konzepts Bogen siehe Definition \ref{bogenkette}).
\end{lem}

\begin{bew}
Sei $b(t) \mathrel{\mathop:}= 1 + (i - 1) t$. Dann benutzen wir den Bogen
$$\left( \tilde X(t), \quad \tilde Y_n(t), \quad \tilde Z_n(t) = 1 + \frac\kappa\mu \left( \tilde Y_n(t) - 1 \right) \right),$$
wobei
$$\tilde X(t) \mathrel{\mathop:}= \begin{cases} b(t)^{-2} \begin{pmatrix} (1 - 6t^2 + 4t^3) & t (1 - t) (2 + 4t) \\ t (1 - t) (6 - 4t) & -(1 - 6t^2 + 4t^3) \end{pmatrix}, &\text{ falls } \lambda = 0 \\ \begin{pmatrix} (1 - 2t) b(t)^{-2} & \lambda \\ \frac4\lambda t (1 - t) b(t)^{-4} & -(1 - 2t) b(t)^{-2} \end{pmatrix}, &\text{ falls } \lambda \neq 0 \end{cases},$$
$$\tilde Y_1(t) \mathrel{\mathop:}= \begin{pmatrix} b(t) & \mu \\ 0 & i b(t) \end{pmatrix}, \quad \tilde Y_2(t) \mathrel{\mathop:}= \begin{pmatrix} b(t) & \mu \\ 0 & -b(t) \end{pmatrix}.$$

Man kann nachprüfen, dass $\tilde X(t)$, $\tilde Y_n(t)$ und $\tilde Z_n(t)$ modulo $\m_K$ sich zu $\begin{pmatrix} 1 & \bar\lambda \\ 0 & 1 \end{pmatrix}$, $\begin{pmatrix} 1 & \bar\mu \\ 0 & 1 \end{pmatrix}$ und $\begin{pmatrix} 1 & \bar\kappa \\ 0 & 1 \end{pmatrix}$ reduzieren. Durch Ausrechnen erhalten wir
$$\tilde X(t)^2 = b(t)^{-4}, \tilde Y_n(t)^4 = b(t)^4 \text{ und } [\tilde Y_n(t), \tilde Z_n(t)] = [\tilde Y_n(t), 1 + \frac\kappa\mu (\tilde Y_n(t) - 1)] = 1,$$
woraus
$$\tilde X(t)^2 \tilde Y_n(t)^4 [\tilde Y_n(t), \tilde Z_n(t)] = 1 \text{ für } 1 \leq n \leq 2$$
folgt. Ferner haben wir $\delta_n = \det \tilde X(t) (\det \tilde Y_n(t))^2 = -b(t)^{-4} (\det \tilde Y_n(t))^2$, also $\delta_1 = -b(t)^{-4} (i b(t)^2)^2 = 1$ und $\delta_2 = -b(t)^{-4} (-b(t)^2)^2 = -1$. Somit gehört der Bogen für $n = 1$ zu $\Spec S^-[1/2]$ und für $n = 2$ zu $\Spec S^+[1/2]$.

Außerdem haben wir $b(0) = 1$ und $b(1) = i$, somit folgen $\tilde X(0) = \tilde X(1) = \tilde X_\lambda$, $\tilde Y_n(0) = \tilde Y_{\mu, n}$ und $\tilde Y_n(1) = \tilde Y_{\mu, n + 2}$, ferner
$$\tilde Z_n(0) = 1 + \frac\kappa\mu (\tilde Y_n(0) - 1) = 1 + \frac\kappa\mu (\tilde Y_{\mu, n} - 1) = \tilde Z_{\mu\kappa, n}$$
und
$$\tilde Z_n(1) = 1 + \frac\kappa\mu (\tilde Y_n(1) - 1) = 1 + \frac\kappa\mu (\tilde Y_{\mu, n + 2} - 1) = \tilde Z_{\mu\kappa, n + 2}.$$
Damit folgt die Behauptung.
\end{bew}

\begin{satz}
\label{spmfrac12intber}
$\Spec S[1/2]$ besteht aus zwei disjunkten irreduziblen Komponenten, nämlich $\Spec S^+[1/2]$ und $\Spec S^-[1/2]$. Insbesondere sind die Ringe $S^+[1/2]$ und $S^-[1/2]$ Integritätsbereiche.
\end{satz}

\begin{bew}
Da nach Lemma \ref{normal} die irreduziblen Komponenten von $\Spec S[1/2]$ disjunkt sind, folgt die erste Behauptung für $\mu = 0$ mit Korollar \ref{punkte1} und Lemma \ref{bogen1}, und für $\mu \neq 0$ mit Korollar \ref{punkte2} und Lemma \ref{bogen2}. Die zweite Behauptung folgt mit der Tatsache, dass die Ringe $S^+[1/2]$ und $S^-[1/2]$ irreduzibel und nach Lemma \ref{normal} als (endliches) direktes Produkt von Integritätsbereichen darstellbar sind.
\end{bew}


\section{Vermutung von Böckle-Juschka}

In Satz \ref{bijektion} dieses Kapitels beantworten wir eine Frage von Böckle-Juschka.


Sei $\1$ der eindimensionale $k$-Vektorraum, auf dem $G_{\Q_2}$ trivial operiert, und sei $D_\1$ der Deformationsfunktor von $\1$. Da $\End_{G_{\Q_2}}(\mathbf{1}) = k$, ist dieser Funktor durch eine vollständige lokale noethersche $\Oc$-Algebra $R_\1$ darstellbar. Wir werden diesen Ring explizit beschreiben. Sei $\psi^{\univ} : G_{\Q_2} \rightarrow R_\1^\times$ die universelle Deformation.

Sei $\Q_2^{\ab}$ der kleinste Unterkörper von $\overline{\Q_2}$, der alle endlichen abelschen Erweiterungen $K$ von $\Q_2$, so dass $[K : \Q_2]$ eine Zweierpotenz ist, enthält. Dann ist $\Gal(\Q_2^{\ab} / \Q_2)$ isomorph zum maximalen pro-2 abelschen Quotient von $G_{\Q_2}$, welchen wir mit $G_{\Q_2}^{\ab}(2)$ bezeichnen. Es folgt aus der lokalen Klassenkörpertheorie, dass $\Q_2^{\ab}(2)$ der kleinste Körper ist, der sowohl die $2$-adische zyklotomische Erweiterung $\Q_2(\mu_{2^\infty})$, als auch die maximale unverzweigte Erweiterung $\Q_2^{\nr}(2)$ in $\Q_2^{\ab}(2)$ enthält. Da $\Q_2(\mu_{2^\infty}) \cap \Q_2^{\nr}(2) = \Q_2$, und $G_{\Q_2}^{\ab}(2)$ abelsch ist, haben wir
\begin{equation}
\label{gq2ab2isom1}
G_{\Q_2}^{\ab}(2) \cong \Gal(\Q_2(\mu_{2^\infty}) / \Q_2) \times \Gal(\Q_2^{\nr}(2) / \Q_2).
\end{equation}
Lokale Klassenkörpertheorie und (\ref{gq2ab2isom1}) führen zu einem Isomorphismus
\begin{equation}
\label{gq2ab2isom2}
G_{\Q_2}^{\ab}(2) \cong \Z_2^\times \times \Z_2 \cong 1 + 4\Z_2 \times \{\pm1\} \times \Z_2.
\end{equation}
Somit können wir Gruppenelemente $\alpha, \beta, \gamma \in G_{\Q_2}$ wählen, so dass deren Bilder in $1 + 4\Z_2 \times \{\pm1\} \times \Z_2$ unter (\ref{gq2ab2isom2}) die Elemente $(5, 1, 0)$, $(1, -1, 0)$ bzw. $(1, 1, 1)$ sind. Da $1 + 4\Z_2$ durch $5$ und $\Z_2$ durch $1$ topologisch erzeugt wird, folgt aus (\ref{gq2ab2isom2}), dass die Bilder von $\alpha$, $\beta$ und $\gamma$ die topologische Gruppe $G_{\Q_2}^{\ab}(2)$ erzeugen.

\begin{prop}
\label{r1isom}
$R_\1 \cong \Oc[[x, y, z]] / ((1 + y)^2 - 1)$.
\end{prop}
\begin{bew}
Sei $(A, \m_A)$ eine lokale artinsche $\Oc$-Algebra mit Restklassenkörper $k$. Dann steht $D_\1(A)$ in einer Bijektion zur Menge der stetigen Gruppenhomomorphismen $\psi : G_{\Q_2} \rightarrow 1 + \m_A$. Da $1 + \m_A$ eine abelsche $2$-Gruppe ist, faktorisiert jeder solche Homomorphismus durch $\psi : G_{\Q_2}^{\ab}(2) \rightarrow 1 + \m_A$. Damit folgt aus (\ref{gq2ab2isom2}), dass die Abbildung $\psi \mapsto (\psi(\alpha) - 1, \psi(\beta) - 1, \psi(\gamma) - 1)$ eine Bijektion zwischen der Menge von solchen $\psi$ und der Menge von Tripeln $(a, b, c) \in \m_a^3$, so dass $(1 + b)^2 = 1$ gilt, induziert, und letztere Menge steht wiederum in einer Bijektion zur Menge der $\Oc$-Algebrenhomomorphismen von $\Oc[[x, y, z]] / ((1 + y)^2 - 1)$ nach $A$.
\end{bew}

\begin{kor}
$R_\1$ ist $\Oc$-torsionsfrei und hat zwei irreduzible Komponenten.
\end{kor}
\begin{bew}
Die erste Behauptung folgt aus der Tatsache, dass $(1 + y)^2 - 1$ in $\Oc[[x, y, z]]$ nicht durch $\varpi$ teilbar ist. Die beiden Komponenten sind gegeben durch $y = 0$ und $y = -2$.
\end{bew}

\begin{defi}
Die Abbildung von einer gerahmten Deformation zu ihrer Determinante induziert eine natürliche Transformation $D^\square \rightarrow D_\1$, und damit einen Homomorphismus von $\Oc$-Algebren $d : R_\1 \rightarrow R^\square$.
\end{defi}

\begin{satz}
\label{bijektion}
Die Abbildung $d : R_\1 \rightarrow R^\square$ induziert eine Bijektion zwischen den irreduziblen Komponenten von $\Spec R^\square$ und $\Spec R_\1$.
\end{satz}
\begin{bew}
Für einen gegebenen Gruppenhomomorphismus $\psi : G_{\Q_2}(2) \rightarrow 1 + \m_K$, wobei $K$ eine endliche Erweiterung von $L$ ist, betrachten wir den Ring
$$S_K[1/2] / (\x_{11} - \psi(x), x_{12}, x_{21}, x_{22}, \y_{11} - \psi(y), y_{21}, y_{22}, z_{12}, z_{21}, \z_{22} - \psi(z) \z_{11}^{-1}),$$
wobei $S_K$ wie $S$ definiert ist, jedoch mit $\Oc_K$ statt $\Oc$. Die Relation $x^2 y^4 [y, z] = 1$ erzwingt, dass wir $\psi(x)^2 \psi(y)^4 = 1$ haben. Durch Ausmultiplizieren erhalten wir dann
\begin{displaymath}
\begin{split}
\tilde X^2 \tilde Y^4 &= \begin{pmatrix} \psi(x)^2 \psi(y)^4 & \psi(x)^2 \y_{12} (\psi(y) + 1) (\psi(y)^2 + 1) + \lambda (\psi(x) + 1) \\ 0 & 1 \end{pmatrix}, \\
[\tilde Y, \tilde Z] &=  \begin{pmatrix} 1 & ((\psi(y) - 1) \kappa + \y_{12} (\psi(z) \z_{11}^{-1} - \z_{11})) \psi(z)^{-1} \z_{11} \\ 0 & 1 \end{pmatrix}.
\end{split}
\end{displaymath}
Damit ist die Relation $\tilde X^2 \tilde Y^4 [\tilde Y, \tilde Z] = 1$ äquivalent zu
\begin{displaymath}
\begin{split}
\psi(x)^2 \y_{12} (\psi(y) + 1) (\psi(y)^2 + 1) + \lambda (\psi(x) + 1) &+ (\psi(y) - 1) \kappa \psi(z)^{-1} \z_{11} \\ &+ \y_{12} (1 - \psi(z)^{-1} \z_{11}^2) = 0.
\end{split}
\end{displaymath}
Da diese Relation im maximalen Ideal von $\Oc_K[[y_{12}, z_{11}]]$ liegt, und der Koeffizient von $y_{12} z_{11}^2$ gleich $\psi(z)^{-1}$, also eine Einheit in $\Oc_K$ ist, folgt nun, dass der genannte Ring vom Nullring verschieden ist. Jedes maximale Primideal dieses Rings korrespondiert zu einem $K$-Punkt $\tilde X, \tilde Y, \tilde Z$ mit $\det \tilde X = \psi(x)$, $\det \tilde Y = \psi(y)$ und $\det \tilde Z = \psi(z)$. Somit induziert die Abbildung $d$ eine Surjektion von maximalen Spektren:
$$\mSpec R^\square[1/2] \rightarrow \mSpec R_\1[1/2].$$
Da der Ring $R_\1[1/2]$ reduziert und Jacobson ist, schließen wir, dass die Abbildung $d : R_\1[1/2] \rightarrow R^\square[1/2]$ injektiv ist. Sei $e = -y/2 \in R_\1[1/2]$, wo $y$ wie in Proposition \ref{r1isom}. Dann haben wir $e^2 = e$, und weil $d$ injektiv ist, ist $d(e)$ ein nichttriviales idempotentes Element in $R^\square[1/2]$. Da
$$R^\square[1/2] \cong S[1/2] \cong S^+[1/2] \times S^-[1/2],$$
und $S^+[1/2]$ und $S^-[1/2]$ nach Satz \ref{spmfrac12intber} beides Integritätsbereiche sind, können wir schließen, dass $d$ eine Bijektion zwischen den irreduziblen Komponenten von $\Spec R^\square[1/2]$ und $\Spec R_\1[1/2]$ induziert. Da $R_\1$ und $R^\square$ beide $\Oc$-torsionsfrei sind, impliziert dies die Behauptung.
\end{bew}


%
%
%
%
%
\end{document}